\documentclass[12pt]{article}
\usepackage[centertags]{amsmath}
\usepackage{amsfonts}
\usepackage{amssymb}
\usepackage{amsthm}
\usepackage{newlfont}
\usepackage{graphicx}

\newfont{\bb}{msbm10 at 12pt}
\def\r{\hbox{\bb R}}

\def\e{\hbox{\bf E}}
\def\t{\hbox{\bf T}}
\def\n{\hbox{\bf N}}
\def\b{\hbox{\bf B}}

\newtheorem{theorem}{Theorem}[section]
\newtheorem{definition}[theorem]{Definition}
\newtheorem{proposition}[theorem]{Proposition}

\newtheorem{lemma}[theorem]{Lemma}

\begin{document}

\title{Time-like Salkowski and anti-Salkowski curves in Minkowski space $\e_1^3$ }
\author{ Ahmad T. Ali\\Mathematics Department\\
 Faculty of Science, Al-Azhar University\\
 Nasr City, 11448, Cairo, Egypt\\
email: atali71@yahoo.com}

\maketitle
\begin{abstract} Salkowski \cite{salkow}, one century ago, introduced a family of curves with constant curvature but non-constant torsion (Salkowski curves) and a family of curves with constant torsion but non-constant curvature (anti-Salkowski curves) in Euclidean 3-space $\e^3$. In this paper, we adapt definition of such curves to time-like curves in Minkowski 3-space $\e_1^3$. Thereafter, we introduce an explicit parametrization of a time-like Salkowski curves and a time-like Anti-Salkowski curves in Minkowski space $\e_1^3$. Also, we characterize them as space curve with constant curvature or constant torsion and whose normal vector makes a constant angle with a fixed line.

\end{abstract}



\section{Introduction}
The Minkowski 3-space $\e_1^3$ is the Euclidean 3-Space $\e^3$ provided with the standard flat metric given by
$$
\langle,\rangle=-dx_1^2+dx_2^2+dx_3^2,
$$
where $(x_1,x_2,x_3)$ is a rectangular coordinate system of $\e_1^3$. If $u=(u_1,u_2,u_3)$ and $v=(v_1,v_2,v_3)$ are arbitrary vectors in $\e_1^3$, we define the (Lorentzian) vector product of $u$ and $v$ as the following:
$$
u\times v=\Bigg|\begin{array}{ccc}
            i & -j & -k \\
            u_1 & u_2 & u_3 \\
            v_1 & v_2 & v_3
          \end{array}\Bigg|.
$$
An arbitrary vector $v\in\e_1^3$ is said space-like if $\langle v,v\rangle>0$ or $v=0$, time-like if $\langle v,v\rangle<0$, and light-like (or null) if $\langle v,v\rangle =0$ and $v\neq0$. The norm (length) of a vector $v$ is given by $\parallel v\parallel=\sqrt{|\langle v,v\rangle|}$.

Given a regular (smooth) curve $\alpha:I\subset\r\rightarrow\e_1^3$, we say that $\alpha$ is space-like (resp.  time-like,  light-like) if all of its velocity vectors $\alpha'(t)$ are  space-like (resp. time-like, light-like). If $\alpha$ is space-like or time-like we say that $\alpha$ is a non-null curve. In such case, there exists a change of the parameter $t$, namely, $s=s(t)$, such that $\parallel\alpha'(s)\parallel=1$.  We say then that $\alpha$ is parameterized by the arc-length parameter.  If the curve $\alpha$ is light-like, the acceleration vector $\alpha''(t)$ must be space-like for all $t$. Then we change the parameter $t$ by $s=s(t)$ in such way that $\parallel \alpha''(s)\parallel=1$ and we say that $\alpha$ is parameterized by the pseudo arc-length parameter. In any of the above cases, we say that $\alpha$ is a unit speed curve \cite{ali,lopez}.

Given a unit speed curve $\alpha$ in Minkowski space $\e_1^3$ it is possible to define a Frenet frame $\{\t(s),\n(s),\b(s)\}$ associated for each point $s$ \cite{kuhn, walr}. Here $\t$, $\n$ and $\b$ are the tangent, normal and binormal vector field, respectively. The geometry of the curve $\alpha$ can be describe by the differentiation of the Frenet frame, which leads to the corresponding Frenet equations. Although different expressions of the Frenet equations appear depending of the causal character of the Frenet trihedron (see the next sections below), we have the concepts of curvature $\kappa$ and torsion $\tau$ of the curve. With this preparatory introduction, we give the following:

Recall that in Euclidean space $\e^3$ a general helix is a curve where the tangent lines make a constant angle with a fixed direction. Helices are characterized by the fact that the ratio $\tau/\kappa$ is constant along the curve \cite{doca}. A slant helix is a curve where the normal lines make a constant angle with a fixed direction \cite{izum}. Slant helices are characterized by the fact that the function $\dfrac{\kappa^2}{(\kappa^2+\tau^2)^{3/2}}\big(\dfrac{\tau}{\kappa}\big)$ is constant \cite{kula}. Salkowski (resp. anti-Salkowski) curves in Euclidean space $\e^3$ are the first known family of curves with constant curvature (resp. torsion) but non-constant torsion (resp. curvature) with an explicit parametrization \cite{monter, salkow}.

In Minkowski space $\e_1^3$, one defines a general helix, slant helix curves in Minkowski space as a similar way. Ferrandez et. al. \cite{ferr} proved that: helices in Minkowski space $\e_1^3$ are characterized by the constancy of the function $\tau/\kappa$ again. Ali and Lopez \cite{ali} proved that: slant helices in Minkowski space $\e_1^3$ are characterized by the constancy of the function $\dfrac{\kappa^2}{(\varepsilon_1\kappa^2+\varepsilon_2\tau^2)^{3/2}}\big(\dfrac{\tau}{\kappa}\big)$, where $\varepsilon_1, \varepsilon_2\in\{-1,1\}$.

In this paper, we define Salkowski curves and anti-Salkowski curves in Minkowski space $\e_1^3$ as a similar way in Euclidean space $\e^3$.o, we introduce the explicit parametrization of a time-like Salkowski curves and time-like anti-Salkowski curves in Minkowski space $\e_1^3$ and we study some characterizations of such curves.

\section{Time-like Salkowski curves and some characterizations}

We suppose that $\alpha$ is a time-like curve. Then $\t^{\prime}(s)\neq 0$ is a space-like vector independent with $\t(s)$. We define the curvature of $\alpha$ at $s$ as $\kappa(s)=|\t^{\prime}(s)|$. The normal vector $\n(s)$ and the binormal $\b(s)$ are defined as
$$
\n(s)=\dfrac{\t^{\prime}(s)}{\kappa(s)}=\dfrac{\psi''}{|\psi''|},\,\,\,\b(s)=\t(s)\times\n(s),
$$
where the vector $\b(s)$ is unitary and space-like. For each $s$, $\{\t,\n,\b\}$ is an orthonormal base of $\e_1^3$ which is called the Frenet trihedron of $\alpha$. We define the torsion of $\alpha$ at $s$ as:
$$
\tau(s)=\langle\n'(s),\b(s)\rangle.
$$
Then the Frenet formula read
\begin{equation}\label{u1}
 \left[
   \begin{array}{c}
     \t' \\
     \n' \\
     \b' \\
   \end{array}
 \right]=\left[
           \begin{array}{ccc}
             0 & \kappa & 0 \\
             \kappa & 0 & \tau \\
             0 & -\tau & 0 \\
           \end{array}
         \right]\left[
   \begin{array}{c}
     \t \\
     \n \\
     \b \\
   \end{array}
 \right],
 \end{equation}
where $g(\t,\t)=-1, g(\n,\n)=g(\b,\b)=1, g(\t,\n)=g(\n,\b)=g(\b,\t)=0$.

We introduce the explicit parametrization of a time-like Salkowski curves in Minkowski space $\e_1^3$ as the following:

\begin{definition} (Time-like Salkowski curves). \label{df-1} For any $m\in R$ with $m>1$, let us define the space curve
\begin{equation}\label{u2}
\begin{array}{ll}
\gamma_m(t)=\dfrac{n}{4m}\Bigg(&\dfrac{1-n}{1+2n}\cosh[(1+2n)t]-\dfrac{1+n}{1-2n}\cosh[(1-2n)t]+2\cosh[t],\\
&\dfrac{1-n}{1+2n}\sinh[(1+2n)t]-\dfrac{1+n}{1-2n}\sinh[(1-2n)t]+2\sinh[t],\\
&\dfrac{1}{m}\sinh[2nt]\Bigg),
\end{array}
\end{equation}
with $n=\dfrac{m}{\sqrt{m^2-1}}$.
\end{definition}
We will name this curve is a time-like Salkowski curve in Minkowski space $\e_1^3$. The geometric elements of the time-like Salkowski curve $\gamma_m$ are the following:

{\bf (1):} $\langle \gamma'_m,\gamma'_m\rangle=-\dfrac{\sinh^2[nt]}{m^2-1}$, so $\|\gamma'_m\|=\dfrac{\sinh[nt]}{\sqrt{m^2-1}}$

{\bf (2):} The arc-length parameter is $s=\dfrac{\cosh[nt]}{m}$.

{\bf (3):} The curvature $\kappa(t)=1$ and the torsion $\tau(t)=\coth[nt]$.

{\bf (4):} The Frenet's frame is
\begin{equation}\label{u3}
\begin{array}{ll}
\t(t)=&\Big(n\cosh[t]\cosh[nt]-\sinh[t]\sinh[nt],\\
&n\sinh[t]\cosh[nt]-\cosh[t]\sinh[nt],\dfrac{n}{m}\cosh[nt]\Big),\\
\n(t)=&\dfrac{n}{m}\Big(\cosh[t],\sinh[t],m\Big),\\
\b(t)=&\Big(\sinh[t]\cosh[nt]-n\cosh[t]\sinh[nt],\\
&\cosh[t]\cosh[nt]-n\sinh[t]\sinh[nt],-\dfrac{n}{m}\sinh[nt]\Big),\\
\end{array}
\end{equation}
From the expression of the normal vector, see Eqs. (\ref{u3}), one can see that the normal indicatrix, or nortrix, of a time-like Salkowski curve in Minkowski space $\e_1^3$ describes a parallel of the unit sphere. The angle between the normal vector and the vector $(0,0,1)$ is constant and equal to $\phi=\mathrm{arccosh}[n]$. This fact is reminiscent of what happens with another important class of curves, the general helices in Minkowski space $\e_1^3$. Such a condition implies that the tangent indicatrix, or tantrix, describes a parallel in the unit sphere.

\begin{lemma} \label{lm-1} Let $\alpha:I\rightarrow E_1^3$ be a time-like curve parameterized by arc-length with $\kappa=1$. Its normal vectors make a constant hyperbolic angle, $\phi$, with a fixed line in space if and only if $\tau(s)=\pm\dfrac{s}{\sqrt{s^2-\tanh^2[\phi]}}$.
\end{lemma}

{\bf Proof:} $(\Rightarrow)$ Let $\textbf{d}$ be the unitary space-like fixed direction which makes a constant hyperbolic angle $\phi$ with the normal vector $\n$. Therefore
\begin{equation}\label{u4}
\langle\n,\textbf{d}\rangle=\cosh[\phi].
\end{equation}
Differentiating Eq. (\ref{u4}) and using Frenet's formula, we get
\begin{equation}\label{u5}
\langle\t+\tau\b,\textbf{d}\rangle=0.
\end{equation}
Therefore,
\begin{equation}\label{u6}
\langle\t,\textbf{d}\rangle=-\tau\langle\b,\textbf{d}\rangle.
\end{equation}
If we put $\langle\b,\textbf{d}\rangle=b$, we can write
\begin{equation}\label{u7}
\textbf{d}=-\tau\,b\,\t+\cosh[\phi]\n+b\,\b.
\end{equation}
From the unitary of the vector $\textbf{d}$ we get $b=\pm\dfrac{\sinh[\phi]}{\sqrt{\tau^2-1}}$. Therefore, the vector $\textbf{d}$ can be written as
\begin{equation}\label{u8}
\textbf{d}=\mp\dfrac{\tau\,\sinh[\phi]}{\sqrt{\tau^2-1}}\t+\cosh[\phi]\n\pm\dfrac{\sinh[\phi]}{\sqrt{\tau^2-1}}\b.
\end{equation}
If we Differentiate Eq. (\ref{u5}) again, we obtain
\begin{equation}\label{u9}
\langle\dot{\tau}\b+(1-\tau^2)\n,\textbf{d}\rangle=0.
\end{equation}
The Eqs. (\ref{u9}) and (\ref{u8}) lead to the differential equation
\begin{equation}\label{u10}
\pm\tanh[\phi]\dfrac{\dot{\tau}}{(\tau^2-1)^{3/2}}+1=0.
\end{equation}
By integration we get
\begin{equation}\label{u11}
\pm\tanh[\phi]\dfrac{\tau}{\sqrt{\tau^2-1}}+s+c=0.
\end{equation}
where $c$ is an integration constant. The integration constant can be subsumed thanks to a parameter change $s\rightarrow s-c$. Finally, to solve (\ref{u11}) with $\tau$ as unknown we get the desired result.

$(\Leftarrow)$ Suppose that $\tau=\pm\dfrac{s}{\sqrt{s^2-\tanh^2[\phi]}}$ and let us consider the vector
\begin{equation}\label{u12}
\textbf{d}=\cosh[\phi]\Big(-s\,\t+\n\mp\sqrt{s^2-\tanh^2[\phi]}\,\b\Big).
\end{equation}
It is easy to prove the vector $\textbf{d}$ is constant i.e., ($\dot{\textbf{d}}=0$) and $\langle\textbf{d},\n\rangle=\cosh[\phi]$.


Once the intrinsic or natural equations of a curve have been determined, the next step is to integrate Frenet's formula with $\kappa=1$ and
\begin{equation}\label{u13}
\tau=\pm\dfrac{s}{\sqrt{s^2-\tanh^2[\phi]}}=\pm\dfrac{\dfrac{s}{\tanh[\phi]}}{\sqrt{\Big(\dfrac{s}{\tanh[\phi]}\Big)^2-1}}=
\pm\coth\Big[\textmd{arccosh}\big[\dfrac{s}{\tanh[\phi]}\big]\Big].
\end{equation}

\begin{theorem} \label{th-1} The time-like curve in Minkowski space $\e_1^3$ with $\kappa=1$ and such that their normal vectors make a constant angle with a fixed line are, up to rigid movements in space or up to the antipodal map, time-like Salkowski curves (see Definition \ref{df-1}).
\end{theorem}

{\bf Proof:} As has been said after Definition \ref{df-1}, the arc-length parameter of time-like Salkowski curves is $s=\int_0^t\|\gamma'_m(u)\|du=\dfrac{1}{m}\cosh[nt]$. Thereafter, $t=\dfrac{1}{n}\textmd{arccosh}[ms]$. In terms of the arc-length curvature and torsion are then
$$
\kappa(s)=1,\,\,\,\,\,\tau(s)=\coth[\textmd{arccosh}[ms]],
$$
the same intrinsic equations, with $m=\coth[\phi]$ and $n=\dfrac{m}{\sqrt{m^2-1}}=\cosh[\phi]$ (compare with the positive case in Eq. (\ref{u13})), as the ones shown in Lemma \ref{lm-1}.

For the negative case in Eq. (\ref{u13}), let us recall that if a curve $\alpha$ has torsion $\tau_{\alpha}$, then the curve $\beta(t)=-\alpha(t)$ has as torsion $\tau_\beta(t)=-\tau_\alpha(t)$, whereas curvature is preserved.

Therefore, the fundamental theorem of curves in Minkowski space states in our situation that, up to a rigid movements or up to the antipodal map, $p\rightarrow-p$, the curves we are looking for are time-like Minkowski curves.

\section{Time-like anti-Salkowski curves }

As an additional material we will show in this section how to build, from a curve in Minkowski space $\e_1^3$ of constant curvature, another curve of constant torsion.

Let us recall that a curve $\alpha:]a,b[\rightarrow\e_1^3$, is 2-regular at a point $t_0$ if $\alpha'(t_0)\neq0$ and if $\kappa_\alpha(t_0)\neq0$.

\begin{lemma} \label{lm-2} Let $\alpha:I\rightarrow E_1^3$ be a regular curve parameterized by arc-length with curvature $\kappa_\alpha$, torsion $\tau_\alpha$ and Frenet's frame $\{\t_\alpha,\n_\alpha,\b_\alpha\}$. Let us consider the curve $\beta(t)=\int_{0}^{t}\t_\alpha(u)\|\b^{'}_\alpha(u)\|\,du$. Then at a parameter $s_\alpha\in I$ such that $\tau_\alpha(s_\alpha)\neq0$, the curve $\beta$ is 2-regular at $s_\beta$ and
$$
\kappa_\beta=\dfrac{\kappa_\alpha}{\tau_\alpha},\,\,\,\tau_\beta=1,\,\,\,\t_\beta=\t_\alpha,\,\,\,
\n_\beta=\n_\alpha,\,\,\,\b_\beta=\b_\alpha.
$$
\end{lemma}

{\bf Proof:} In order to obtain the tangent vector of $\beta$ let us compute
\begin{equation}\label{u14}
\t_\beta(s_\beta)=\dot{\beta}(s_\beta)=\dfrac{d\beta}{dt}\dfrac{dt}{ds_\beta}=\t_\alpha\|\b'_\alpha(t)\|\dfrac{dt}{ds_\beta}.
\end{equation}
From the above equation, we get
\begin{equation}\label{u15}
\dfrac{ds_\beta}{dt}=\|\b'_\alpha(t)\|=\Big\|\dfrac{\b_\alpha}{ds_\alpha}\dfrac{ds_\alpha}{dt}\Big\|=\tau_\alpha\dfrac{ds_\alpha}{dt},
\end{equation}
and
\begin{equation}\label{u151}
\t_\beta(s_\beta)=\t_\alpha(s_\alpha).
\end{equation}
Differentiation the above equation along with Frenet's Eqs. (\ref{u1}) we obtain
\begin{equation}\label{u16}
\dot{\t}_\beta(s_\beta)=\dfrac{d\t_\alpha}{ds_\alpha}\,\dfrac{ds_\alpha}{dt}\,\dfrac{dt}{ds_\beta}.
\end{equation}
Using Frenet's Eqs. (\ref{u1}) and Eq. (\ref{u15}), the above equation can be written as
\begin{equation}\label{u17}
\kappa_\beta\,\n_\beta(s_\beta)=\dfrac{\kappa_\alpha}{\tau_\alpha}\,\n_\alpha(s_\alpha)
\end{equation}
From the above equation, we get
\begin{equation}\label{u18}
\kappa_\beta=\dfrac{\kappa_\alpha}{\tau_\alpha},
\end{equation}
and
\begin{equation}\label{u181}
\n_\beta(s_\beta)=\n_\alpha(s_\alpha).
\end{equation}
So that, we have
\begin{equation}\label{u19}
\b_\beta(s_\beta)=\t_\beta(s_\beta)\times\n_\beta(s_\beta)=\t_\alpha(s_\alpha)\times\n_\alpha(s_\alpha)=\b_\alpha(s_\alpha).
\end{equation}
Differentiating the above equation with respect to $s_\beta$ we get $\tau_\beta=1$.

Let us apply the previous result to the time-like Salkowski curve $\gamma_m$ defined in Eq. (\ref{u2}) we have the explicit parametrization of a time-like anti-Salkowski curve as the following:
\begin{equation}\label{u111}
\begin{array}{ll}
\beta_m(t)=\dfrac{n}{4m}\Bigg(&\dfrac{n-1}{2n+1}\sinh[(1+2n)t]-\dfrac{n+1}{2n-1}\sinh[(1-2n)t]+2n\sinh[t],\\
&\dfrac{n-1}{2n+1}\cosh[(1+2n)t]-\dfrac{n+1}{2n-1}\cosh[(1-2n)t]+2n\cosh[t],\\
&-\dfrac{1}{m}(\sinh[2nt]+2nt)\Bigg),
\end{array}
\end{equation}
where, as for time-like Salkowski curves, $n=\dfrac{m}{\sqrt{m^2-1}}$ and $m>1$. Let us call these curves by the name time-like anti-Salkowski curves. The presence of the non-trigonometric term $2nt$ in the third component of $\beta_m$ makes that the change of variable studied in Section 2 for time-like Salkowski curves does not work for anti-Salkowski curves.

Applying Lemma \ref{lm-2} we get the following proposition:

\begin{proposition} \label{pr-1} The curve $\beta_m$ in Eq. (\ref{u111}) are curves of constant torsion equal to 1 and non-constant curvature equal to $\tanh[nt]$.
\end{proposition}

Finally, we stated here the following Lemma:

\begin{lemma} \label{lm-3} Let $\alpha:I\rightarrow E_1^3$ be a regular curve parameterized by arc-length with curvature $\kappa_\alpha$, torsion $\tau_\alpha$ and Frenet's frame $\{\t_\alpha,\n_\alpha,\b_\alpha\}$. Let us consider the curve $\beta(t)=\int_{0}^{t}\t_\alpha(u)\|\t^{\,\prime}_\alpha(u)\|\,du$. Then at a parameter $s_\alpha\in I$ such that $\kappa_\alpha(s_\alpha)\neq0$, the curve $\beta$ is 2-regular at $s_\beta$ and
$$
\kappa_\beta=1,\,\,\,\tau_\beta=\dfrac{\tau_\alpha}{\kappa_\alpha},\,\,\,\t_\beta=\t_\alpha,\,\,\,
\n_\beta=\n_\alpha,\,\,\,\b_\beta=\b_\alpha.
$$
\end{lemma}

{\bf Proof:} The proof of this Lemma is the same as the proof of Lemma \ref{lm-2}.

\begin{theorem} \label{th-2} The space curve with $\tau=1$ and such that their normal vectors makes a constant angle with a fixed line are the anti-Salkowski curves defined in Eq. (\ref{u111}).
\end{theorem}

{\bf Proof:} Let $\alpha$ be a curve with $\tau=1$ and let $\beta(t)=\int_{0}^{t}\t_\alpha(u)\|\t'_\alpha(u)\|\,du$. By Lemma \ref{lm-3}, $\beta$ is a curve with constant curvature $\kappa=1$, non-constant torsion $\tau=\dfrac{1}{\kappa_\alpha}$ and with the same normal vector. Therefore, $\beta$ is a Salkowski curve and $\alpha$ is an anti-Salkowski curve.


\end{document}